\documentclass{amsart}

\usepackage{mathtools}
\usepackage[bigdelims,vvarbb]{newtxmath} 
\usepackage{comment}

\newtheorem{theorem}{Theorem}[section]

\theoremstyle{definition}

\newtheorem{remark}[theorem]{Remark}

\newcommand{\sameorder}{\asymp}
\newcommand{\dx}{\mathrm{d}}

\newcommand{\M}{\mathcal{M}}
\newcommand{\C}{\mathbb{C}}
\newcommand{\N}{\mathbb{N}}
\newcommand{\R}{\mathbb{R}} 
\newcommand{\Y}{\mathbb{Y}}
\newcommand{\doublesum}{\mathop{\sum\sum}}
\newcommand{\Stilde}{\widetilde{S}}

\newcommand{\Odipg}[2]{\mathcal{O}_{#1}\bigl(#2\bigr)}
\newcommand{\Odipm}[2]{\mathcal{O}_{#1} (#2)}

\allowdisplaybreaks

\renewcommand{\qedsymbol}{$\square$}
\newenvironment{Proof}[1][Proof]{\par\noindent\textbf{#1.}~}
{\hfill\qedsymbol\smallskip\par}
\newtheorem{Theorem}{Theorem}
\newtheorem{Lemma}{Lemma}

\begin{document}


\title[A Ces\`aro Average for an additive problem with prime powers]{A Ces\`aro Average for an additive problem\\ with prime powers} 
\date{\today}

\author{Alessandro Languasco}
\address{Dipartimento di Matematica ``Tullio Levi-Civita'',
Universit\`a di Padova \\
Via Trieste 63, 35121 Padova, Italy \\
E-mail: alessandro.languasco@unipd.it}

\author{Alessandro Zaccagnini}
\address{Dipartimento di Scienze Matematiche, Fisiche e Informatiche,
Universit\`a di Parma \\
Parco Area delle Scienze, 53/a, 43124 Parma, Italy \\
E-mail: alessandro.zaccagnini@unipr.it}


\maketitle

\begin{abstract}
In this paper we extend and improve our results on weighted averages
for the number of representations of an integer as a sum of two powers
of primes, that appeared in \cite{LanguascoZ2015a} (see also Theorem
2.2 of \cite{Languasco2016a}).
Let  $1\le \ell_1 \le \ell_2$ be two integers, $\Lambda$ be the von Mangoldt function and
\(
r_{\ell_1,\ell_2}(n) =   \sum_{m_1^{\ell_1} + m_2^{\ell_2}= n} \Lambda(m_1) \Lambda(m_2)
\)
be the weighted counting function for the number of representation
of an integer as a sum of two prime powers. Let $N \geq 2$
be an integer. We prove that the Ces\`aro average of weight $k > 1$
of $r_{\ell_1,\ell_2}$ over the interval $[1, N]$ has a development as
a sum of terms depending explicitly on the zeros of the Riemann
zeta-function.
\end{abstract}

\section{Introduction} 
We continue our recent work on the number of representations of an
integer as a sum of primes.
In \cite{LanguascoZ2012a} we studied the \emph{average} number of
representations of an integer as a sum of two primes, whereas in
\cite{LanguascoZ2012b} we considered individual integers. 
In \cite{LanguascoZ2015a},
see also Theorem 2.2 of \cite{Languasco2016a}, we  studied a Ces\`aro
weighted partial \emph{explicit} formula for Goldbach numbers.
Here we generalise and improve such last result
by working  on the Ces\`aro weighted counting function for the number of representation
of an integer as a sum of two prime powers.
We let $1\le \ell_1 \le \ell_2$ be two integers and set
\[
  r_{\ell_1,\ell_2}(n)
  =
  \sum_{m_1^{\ell_1} + m_2^{\ell_2}= n} \Lambda(m_1) \Lambda(m_2).
\]
We also use the following convenient abbreviations for the various
terms of the development:
\begin{align}
\notag
  \M_{1, k, \ell_1, \ell_2}(N)
  &=
  \frac{N^{ 1/\ell_1 + 1/\ell_2}}{\Gamma(k +1+ 1/\ell_1 + 1/\ell_2 )}
  \frac{\Gamma(1/\ell_1)\Gamma(1/\ell_2)}{\ell_1\ell_2}, \\
  \M_{2, k}(N)
\notag
  &=
  \frac{\log^2(2\pi)}{\Gamma(k + 1)}, \\
\label{def-M3}
  \M_{3, k, \ell}(N)
  &=
  - \log(2\pi)
  \frac{N^{1/\ell}}{\Gamma(k +1+1/\ell)}
  \frac{ \Gamma(1/\ell)}{\ell}
  +
  \log(2\pi)
  \sum_{\rho}  \Gamma\Bigr(\frac{\rho}{\ell}\Bigl)  
  \frac{N^{\rho/\ell}}{\Gamma( k + 1+  \rho/\ell)}, \\
\label{def-M4}
  \M_{4, k, \ell_1, \ell_2}(N)
  &=
  - N^{ 1/\ell_2 }\frac{ \Gamma(1/\ell_2)}{\ell_2}
  \sum_{\rho}  \Gamma\Bigr(\frac{\rho}{\ell_1}\Bigl)  
  \frac{N^{\rho/\ell_1 }}{\Gamma( k + 1+1/\ell_2 + \rho/\ell_1)}, \\
\label{def-M5}
  \M_{5, k, \ell_1, \ell_2}(N)
  &=
  \sum_{\rho_1} \sum_{\rho_2}
    \frac{ \Gamma ( \rho_1/\ell_1) \Gamma ( \rho_2/\ell_2)}
  {\Gamma(k + 1+ \rho_1/\ell_1 +\rho_2/\ell_2)}
     N^{ \rho_1/\ell_1 + \rho_2/\ell_2}.
\end{align}
Here $\rho$, with or without subscripts, runs over the non-trivial
zeros of the Riemann zeta-function $\zeta$ and $\Gamma$ is Euler's
function.
The main result of the paper is the following theorem.

\begin{Theorem}
\label{Cesaro-gen-average}
Let $1\le \ell_1 \le \ell_2$ be two integers, and 
$N$ be a positive integer.
For $k > 1$ we have
\begin{align*}
  \sum_{n \le N} r_{\ell_1,\ell_2}(n) \frac{(1 - n/N)^k}{\Gamma(k + 1)}
  &=
  \M_{1, k, \ell_1, \ell_2}(N) + \M_{2, k}(N) + \M_{3, k, \ell_1}(N) + \M_{3, k, \ell_2}(N) \\
  &\qquad+
  \M_{4, k, \ell_1, \ell_2}(N) + \M_{4, k, \ell_2, \ell_1}(N) + \M_{5, k, \ell_1, \ell_2}(N) \\
  &\qquad+
  \Odipm{k,\ell_1,\ell_2}{N^{- 1/2+1/\ell_1}}.
\end{align*}
\end{Theorem}

Clearly, depending on the size of $\ell_1,\ell_2$, 
some of the previous listed terms should be included
in the error term.
We remark that the double series over zeros in \eqref{def-M5} 
converges absolutely for $k > 1 / 2$, and it seems reasonable to
believe that the stated equality holds for the same values of $k$,
possibly with a weaker error term, although the bound $k>1$ appears in
several places of the proof and it seems to be the limit of the
method.

Theorem \ref{Cesaro-gen-average}  generalises and improves our Theorem in \cite{LanguascoZ2015a},
see also Theorem 2.2 of \cite{Languasco2016a}, which corresponds to the case $\ell_1=\ell_2=1$.
In fact in this case Theorem \ref{Cesaro-gen-average} leads to 
\begin{align*}
\notag
  &\sum_{n \le N} r_{G}(n) \frac{(1 - n/N)^k}{\Gamma(k + 1)}
  =
   \frac{N^{2}}{\Gamma(k +3)}  
   -  2 
  \sum_{\rho} 
  \frac{ \Gamma (  \rho)  }{\Gamma( k + 2+ \rho )}N^{\rho +1}
   - 2\log(2\pi)
   \frac{N }{\Gamma(k +2 )} 
   \\&
    + 2  \log(2\pi)
  \sum_{\rho}   
  \frac{\Gamma (  \rho )}{\Gamma( k + 1+  \rho )}N^{  \rho }
+
     \sum_{\rho_1} \sum_{\rho_2}
    \frac{ \Gamma ( \rho_1 ) \Gamma ( \rho_2 )}
  {\Gamma(k + 1 + \rho_1  +\rho_2)}
     N^{ \rho_1  + \rho_2 }
  +
   \Odipm{k}{N^{1/2}},
\end{align*} 
where \(
r_{G}(n) =
r_{1,1}(n) =   \sum_{m_1 + m_2 = n} \Lambda(m_1) \Lambda(m_2),
\)
that is, we are now able to detect the term $\M_{3,k,1}$.
Very recently Br\"udern, Kaczorowski and Perelli \cite{BrudernKP2017}
found the full explicit formula for $k > 0$.
We point out that Theorem \ref{Cesaro-gen-average}  covers other interesting and classical cases  like the 
sum of two prime squares ($\ell_1=\ell_2=2$) or a prime and a prime square
($\ell_1=1,\ell_2=2$).

We recall that our method  is based on a formula due to Laplace
\cite{Laplace1812}, namely
\begin{equation}
\label{Laplace-transf}
  \frac 1{2 \pi i}
  \int_{(a)} v^{-s} e^v \, \dx v
  =
  \frac1{\Gamma(s)},
\end{equation}
where $\Re(s) > 0$ and $a > 0$, see, e.g., formula 5.4(1) on page 238
of \cite{ErdelyiMOT1954a}.
We will need the general case of \eqref{Laplace-transf},
which can be found in de Azevedo Pribitkin \cite{Azevedo2002},
formulae (8) and (9):
\begin{equation}
\label{Laplace-eq-1}
  \frac1{2 \pi}
  \int_{\R} \frac{e^{i D u}}{(a + i u)^s} \, \dx u
  =
  \begin{cases}
    \dfrac{D^{s - 1} e^{- a D}}{\Gamma(s)}
    & \text{if $D > 0$,} \\
    0
    & \text{if $D < 0$,}
  \end{cases}
\end{equation}
which is valid for $\sigma = \Re(s) > 0$ and $a \in \C$ with
$\Re(a) > 0$, and
\begin{equation}
\label{Laplace-eq-2}
  \frac1{2 \pi}
  \int_{\R} \frac 1{(a + i u)^s} \, \dx u
  =
  \begin{cases}
    0     & \text{if $\Re(s) > 1$,} \\
    1 / 2 & \text{if $s = 1$,}
  \end{cases}
\end{equation}
for $a \in \C$ with $\Re(a) > 0$.
Formulae \eqref{Laplace-eq-1}-\eqref{Laplace-eq-2} enable us to write
averages of arithmetical functions by means of line integrals as we
will see in \S\ref{settings} below.

The improvement we get in Theorem \ref{Cesaro-gen-average} 
follows using Lemma \ref{Linnik-lemma-gen}
below which is a generalised and refined
version of Lemma 4.1 of
 \cite{LanguascoZ2015a},
see also Lemma 5.1 of \cite{Languasco2016a}. 
In fact Lemma \ref{Linnik-lemma-gen}
can be also used to generalise and improve 
our result in \cite{LanguascoZ2013a}
about the Hardy-Littlewood numbers
to the $p^\ell+m^{2}$, $\ell\ge 1$, problem;  
we will discuss this case in \cite{LanguascoZ2017d}.

\section{Settings}
\label{settings}

Let $\ell \ge 1$, $1\le \ell_1 \le \ell_2$  be integral numbers  and
\begin{equation}
\label{Stilde-def}
  \Stilde_\ell(z)
  =
  \sum_{m \ge 1}\Lambda(m) e^{- m^\ell z},
\end{equation}
where $z = a + i y$ with $y \in \R$ and real $a > 0$.
Moreover let us define the density of the problem as
\begin{equation}
\label{lambda-def}
\lambda =  1/\ell_1+1/\ell_2.
\end{equation}

We recall that the Prime Number Theorem (PNT) is equivalent 
to the statement
\begin{equation}
\label{PNT-equiv}
  \Stilde_\ell(a)
  \sim
  \frac{\Gamma(1/\ell)}{\ell a^{1/\ell}}
  \qquad\text{for $a \to 0+$.}
\end{equation}
By \eqref{Stilde-def} we have
\[
  \Stilde_{\ell_1}(z) \Stilde_{\ell_2}(z) 
  =
  \sum_{n \ge 1} r_{\ell_1,\ell_2}(n) e^{-n z}.
\]
Hence, for $N \in \N$ with $N > 0$ and $a > 0$ we have
\begin{equation}
\label{starting-point}
 \frac 1{2 \pi i}
  \int_{(a)} e^{N z} z^{- k - 1} \Stilde_{\ell_1}(z) \Stilde_{\ell_2}(z) \, \dx z
  =
  \frac 1{2 \pi i}
  \int_{(a)} e^{N z} z^{- k - 1}
    \sum_{n \ge 1} r_{\ell_1,\ell_2}(n) e^{- n z} \, \dx z.
\end{equation}
Since
\[
  \sum_{n \ge 1} \vert r_{\ell_1,\ell_2}(n) e^{- n z} \vert 
  =
 \Stilde_{\ell_1}(a) \Stilde_{\ell_2}(a)
  \sameorder_{\ell_1,\ell_2}
  a^{-\lambda}
\]
by \eqref{PNT-equiv}, where $f\sameorder g$ means $g \ll f \ll g$, we
can exchange the series and the line integral in
\eqref{starting-point} provided that $k>0$.
In fact, if $z = a + i y$, taking into account the estimate
\begin{equation}
\label{z^-1}
  \vert z \vert^{-1}
  \sameorder
  \begin{cases}
    a^{-1}   &\text{if $\vert y \vert \le a$,} \\
    \vert y \vert^{-1} &\text{if $\vert y \vert \ge a$,}
  \end{cases}
\end{equation}
we have
\[
  \vert e^{N z} z^{- k - 1}\vert  
  \sameorder
   e^{N a}
  \begin{cases}
    a^{- k - 1} &\text{if $\vert y \vert \le a$,} \\
    \vert y \vert^{- k - 1} &\text{if $\vert y \vert \ge a$,}
  \end{cases}
\]
and hence, recalling \eqref{PNT-equiv}, we obtain
%
%
%
\[
  \int_{(a)} \vert e^{N z} z^{- k - 1} \vert \,
    \Bigl\vert
      \sum_{n \ge 1} r_{\ell_1,\ell_2}(n) e^{- n z}
    \Bigr\vert \, \vert \dx z \vert
  \ll
  a^{-\lambda} e^{N a}
  \Bigl[
    \int_0^a a^{- k - 1} \, \dx y
    +
    \int_a^{+\infty} y^{- k - 1} \, \dx y
  \Bigl],
\]
which is $ \ll_k a^{-\lambda-k} e^{N a}$, but the rightmost integral
converges only for $k > 0$.
Using \eqref{Laplace-eq-1} for $n \ne N$ and \eqref{Laplace-eq-2} for
$n = N$, we see that for $k > 0$ the right-hand side of
\eqref{starting-point} is
\begin{align*}
  &=
  \sum_{n \ge 1} r_{\ell_1,\ell_2}(n)
    \Bigl(
      \frac 1{2 \pi i}
      \int_{(a)} e^{(N - n) z} z^{- k - 1} \, \dx z
    \Bigr)  
  =
  \sum_{n \le N} r_{\ell_1,\ell_2}(n) \frac{(N - n)^k}{\Gamma(k + 1)}.
\end{align*}

\begin{remark}
As in \cite{LanguascoZ2015a}  the previous computation reveals
that we can not get rid of the Ces\`aro weight in our method since,
for $k = 0$, it is not clear whether the integral on the right hand side
of \eqref{starting-point} converges absolutely or not.
\end{remark}

Summing up, for $a>0$ and $k > 0$ we have
\begin{equation*} 
  \sum_{n \le N}
    r_{\ell_1,\ell_2}(n) \frac{(N - n)^k}{\Gamma(k + 1)}
  =
 \frac 1{2 \pi i}
  \int_{(a)} e^{N z} z^{- k - 1} \Stilde_{\ell_1}(z) \Stilde_{\ell_2}(z)   \, \dx z,
\end{equation*}
where $N \in \N$ with $N > 0$.
This is the fundamental relation for the method.

\section{Inserting zeros}

In this section we need $k > 1$.
By Lemma \ref{Linnik-lemma-gen} below we have
\[
  \Stilde_{\ell}(z)  
  =
 \frac{\Gamma(1/\ell)}{\ell z^{1/\ell}}
- 
\frac{1}{\ell}\sum_{\rho}z^{-\rho/\ell}\Gamma\Bigl(\frac{\rho}{\ell}\Bigr) 
- \log (2\pi)
  +
  E(a,y,\ell)=
  M(\ell,z) +E(a,y,\ell),
\]
say, where $E(a,y,\ell)$ satisfies \eqref{expl-form-err-term-strong}.
Hence
\begin{align*}
  \Stilde_{\ell_1}(z)\Stilde_{\ell_2}(z)
 & =
 M(\ell_1,z) M(\ell_2,z) 
  +
  E(a,y,\ell_1) E(a,y,\ell_2) 
  \\&
  \qquad+
  E(a,y,\ell_2)  M(\ell_1,z) +    E(a,y,\ell_1) M(\ell_2,z).
\end{align*}
We have
$\vert  M(\ell,z)  \vert
   =
   \bigl\vert \Stilde_\ell(z) - E(a,y,\ell) \bigr\vert
  \le
  \Stilde_\ell(a) + \bigl\vert E(a,y,\ell) \bigr\vert
  \ll_\ell
  a^{-1/\ell}
  +
  \bigl\vert E(a,y,\ell) \bigr\vert$
by \eqref{PNT-equiv} again, so that
\begin{align}
\notag
 \Stilde_{\ell_1}(z) \Stilde_{\ell_2}(z)
 &=
 M(\ell_1,z) M(\ell_2,z) 
 +
  \Odipg{\ell_1,\ell_2}{\bigl\vert  E(a,y,\ell_1) E(a,y,\ell_2) \bigr\vert}
 \\&
 \label{squaring-out} 
  \qquad+
  \Odipg{\ell_1,\ell_2}{\bigl\vert  E(a,y,\ell_2) \bigr\vert a^{-1/\ell_1}
    +
  \bigl\vert  E(a,y,\ell_1) \bigr\vert a^{-1/\ell_2}}\ ,
\end{align}
choosing $0<a \le 1$,  since $1\le\ell_1\le \ell_2$.
Recalling \eqref{z^-1} and taking into account
\eqref{expl-form-err-term-strong}, for $k>1$ we have
\begin{align*}
  \int_{(a)} &\bigl\vert  E(a,y,\ell_1) E(a,y,\ell_2) \bigr\vert 
    \vert e^{N z} \vert \, \vert z \vert^{- k - 1} \, \vert \dx z \vert 
  \\&
  \ll_{\ell_1,\ell_2}
  e^{N a}
  \int_0^a a^{- k} \, \dx y
  +
  e^{N a}
  \int_a^{+\infty} y^{- k}(1+ \log^2(y/a))^2 \, \dx y \\
  &\ll_{k,\ell_1,\ell_2}
  e^{N a} a^{- k + 1} +  e^{N a} a^{- k + 1}
  \int_1^{+\infty} v^{- k} (1 + \log^2v)^2 \, \dx v 
  \ll_{k,\ell_1,\ell_2}
  e^{N a} a^{- k + 1}.
\end{align*}
Choosing $a = 1 / N$, the error term is $\ll_{k,\ell_1,\ell_2} N^{k - 1}$ for $k > 1$.
For $a = 1 / N$, by \eqref{z^-1} and \eqref{expl-form-err-term-strong},
the second remainder term in \eqref{squaring-out}  for $k > 1/2$ is
\begin{align*}
  &\ll_{\ell_1,\ell_2}
 N^{1/\ell_1}
  \int_{(\frac{1}{N})} \vert E(y,1/N, \ell_2) \vert\vert e^{N z}\vert
    \vert z\vert ^{- k - 1} \, \vert \dx z \vert \\
  &\ll_{\ell_1,\ell_2}
   N^{1/\ell_1}
  \int_0^{1/N} N^{ k + 1/2} \, \dx y
  +
   N^{1/\ell_1}
  \int_{1/N}^{+\infty} y^{- k - 1/2} \log^2(Ny) \, \dx y \\
  &\ll_{k,\ell_1,\ell_2}
  N^{k - 1/2+1/\ell_1} 
  +
  N^{k - 1/2+1/\ell_1} 
  \int_{1}^{+\infty}  v^{- k - 1/2} \log^2 v \, \dx v 
  \ll_{k,\ell_1,\ell_2}
  N^{k - 1/2+1/\ell_1} .
\end{align*}
Arguing analogously, it is easy to see that the remaining term is 
$\ll_{k,\ell_1,\ell_2} N^{k - 1/2+1/\ell_2} $.
  
With a little effort we can give an explicit dependence on $k$ for the
implicit constants in the last three estimates.

Hence, by  \eqref{lambda-def} and \eqref{starting-point} we have
\begin{align*}
  \sum_{n \le N} r_{\ell_1,\ell_2}(n) \frac{(N - n)^k}{\Gamma(k + 1)}
  &=
  \frac 1{2 \pi i}
  \int_{(\frac{1}{N})}
    e^{N z} z^{- k - 1}
    M(\ell_1,z) M(\ell_2,z) 
    \, \dx z
  +
  \Odipm{k,\ell_1,\ell_2}{N^{k - 1/2+1/\ell_1}} \\
  &=
  I_1(N; \ell_1, \ell_2, k)
  +
  I_2(N; k)
  +
  I_3(N; \ell_1, k) + I_3(N; \ell_2, k) \\
  &\qquad+
  I_4(N; \ell_1, \ell_2, k)
  +
  I_4(N; \ell_2, \ell_1, k)
  +
  I_5(N; \ell_1, \ell_2, k) \\
  &\qquad+
  \Odipm{k,\ell_1,\ell_2}{N^{k - 1/2+1/\ell_1}},
\end{align*}
say, where
\begin{align*}
  I_1(N; \ell_1, \ell_2, k)
  &=
  \frac 1{2 \pi i} \frac{\Gamma(1/\ell_1)\Gamma(1/\ell_2)}{\ell_1\ell_2}
  \int_{(\frac{1}{N})} e^{N z} z^{- k - 1-\lambda} \, \dx z, \\
  I_2(N; k)
  &=
  \frac {\log^2(2\pi)}{2 \pi i} 
    \int_{(\frac{1}{N})} e^{N z} z^{- k - 1} \, \dx z, \\
  I_3(N; \ell, k)
  &=
  \frac {\log(2\pi)}{2 \pi i}
  \Bigl\{
    -\frac{ \Gamma(1/\ell)}{ \ell}
    \int_{(\frac{1}{N})} e^{N z} z^{- k - 1- 1/\ell} \, \dx z
  +
  \int_{(\frac{1}{N})} e^{N z} z^{- k - 1 }
    \sum_{\rho} z^{-\frac{\rho}{\ell}} \Gamma\Bigr(\frac{\rho}{\ell}\Bigl) \, \dx z
  \Bigr\}, \\
  I_4(N; \ell_1, \ell_2, k)
  &=
  -
  \frac 1{2\pi i}\frac{ \Gamma(1/\ell_1)}{ \ell_1}
  \int_{(\frac{1}{N})} e^{N z} z^{- k - 1 -1/\ell_1}
    \sum_{\rho} z^{-\frac{\rho}{\ell_2}} \Gamma\Bigr(\frac{\rho}{\ell_2}\Bigl) \, \dx z,
\\ 
  I_5(N; \ell_1, \ell_2, k)
  &=
  \frac 1{2 \pi i}
  \int_{(\frac{1}{N})}
    e^{N z} z^{- k - 1}
    \sum_{\rho_1} \sum_{\rho_2}
      z^{-\frac{\rho}{\ell_1}-\frac{\rho}{\ell_2}}\Gamma\Bigr(\frac{\rho_1}{\ell_1}\Bigl) \Gamma\Bigr(\frac{\rho_2}{\ell_2}\Bigl)
    \, \dx z.
\end{align*}
The evaluation of the integrals $I_j$ is a straightforward application
of \eqref{Laplace-transf} with $s = N z$, except that the interchange
of the series with the integrals needs to be justified: see
\S\ref{first-exchange}-\ref{exchange-double-sum-rhos} for a proof that
this is in fact permitted when $k > 1$.
The proof that the double sum over zeros converges absolutely for
$k > 1 / 2$ is given in \S\ref{sec:double-sum} below.
Combining the resulting expressions and dividing through by $N^k$
we get Theorem~\ref{Cesaro-gen-average}.

\section{Lemmas}
We recall some basic facts in complex analysis.
First, if $z = a + i y$ with $a > 0$, we see that for complex $w$ we
have
\begin{align*}
  z^{-w}
  &=
  \vert z \vert^{-w} \exp( - i w \arctan(y / a))
  =
  \vert z \vert^{-\Re(w) - i \Im(w)} \exp( (- i \Re(w) + \Im(w)) \arctan(y / a))
\end{align*}
so that
\begin{equation}
\label{z^w}
  \vert z^{-w} \vert
  =
  \vert z \vert^{-\Re(w)} \exp(\Im(w) \arctan(y / a)).
\end{equation}
We also recall that, uniformly for $x \in [x_1, x_2]$, with $x_1$ and
$x_2$ fixed, and for $|y| \to +\infty$, by the Stirling formula
(see, e.g., Titchmarsh \cite[\S4.42]{Titchmarsh1988}) we have
\begin{equation}
\label{Stirling}
  \vert \Gamma(x + i y) \vert
  \sim
  \sqrt{2 \pi}
  e^{- \pi |y| / 2} |y|^{x - 1 / 2}.
\end{equation}

The following lemma generalizes and improves Lemma 4.1 of
 \cite{LanguascoZ2015a},
see also Lemma 5.1 of \cite{Languasco2016a}.
The improvement depends on the fact that the constant term
$\log(2\pi)$ is now explicit since we realised that, in the application,
this term leads, in some cases, to a non-trivial contribution 
in the final result.
We follow the
line of the proof in  \cite{LanguascoZ2015a},  but, in some cases,  
the integration path has to be changed; for  clarity
we repeat the whole argument. 
\begin{Lemma} 
\label{Linnik-lemma-gen}
Let $\ell\ge 1$ be an integer,  $z = a + iy$, where $a > 0$ and $y \in \R$.
Then
\begin{equation*}  
\Stilde_{\ell}(z)  
= 
\frac{\Gamma(1/\ell)}{\ell z^{1/\ell}}
- 
\frac{1}{\ell}\sum_{\rho}z^{-\rho/\ell}\Gamma\Bigl(\frac{\rho}{\ell}\Bigr) 
- \log (2\pi)
   +
  E(a,y,\ell),
\end{equation*}
where $\rho = \beta + i\gamma$ runs over the non-trivial zeros of
$\zeta(s)$ and
\begin{equation}
\label{expl-form-err-term-strong}
  E(a,y,\ell)
  \ll_\ell
  \vert z \vert^{1/2}
  \begin{cases}
    1 & \text{if $\vert y \vert \leq a$} \\
    1 +\log^2 (\vert y\vert/a) & \text{if $\vert y \vert > a$.}
  \end{cases}
\end{equation}
\end{Lemma}

\begin{Proof}
Following the line of Hardy and Littlewood, see
\cite[\S2.2]{HardyL1916}, \cite[Lemma 4]{HardyL1923} and of \S4 in
Linnik \cite{Linnik1946}, we have that
\begin{align}
\notag
   \Stilde_\ell(z)
  &=
\frac{\Gamma(1/\ell)}{\ell z^{1/\ell}}
- 
\frac{1}{\ell}\sum_{\rho}z^{-\rho/\ell}\Gamma\Bigl(\frac{\rho}{\ell}\Bigr) 
  -
  \frac{\zeta'}{\zeta}(0)
  -
  \sum_{m=1}^{\ell/4} \Gamma\Bigl(-\frac{2m}{\ell}\Bigr) z^{2m/\ell}
  \\
 \label{Mellin2}
   &\hskip1cm
  -
  \frac{1}{2\pi i}
  \int_{\mathcal{L}_\ell} 
  \frac{\zeta'}{\zeta}(\ell w) \Gamma(w)z^{-w} \, \dx w,
\end{align}
where $\mathcal{L}_\ell$ is the vertical line $\Re(w)=-1/2$ if $4\nmid \ell$
and it is $\{-1/2+it\colon \vert t \vert > C\} \cup \{-1/2+it\colon 1/\ell \le \vert t \vert \le C\} \cup \gamma_{\ell}$
otherwise,
  $C>1/\ell$ is an absolute constant to be chosen later and $\gamma_{\ell}$
is the right half-circle centred in $-1/2$ of radius $1/\ell$.

Now we estimate the integral in \eqref{Mellin2}. 
Assume $4\nmid \ell$. 
Writing $w=-1/2+it$, we have
$\vert (\zeta'/\zeta)(\ell w)\vert \ll_\ell \log (\vert t \vert +2)$,
$\vert z^{-w} \vert = \vert z \vert^{1 / 2} \exp(t \arctan( y / a ))$
by \eqref{z^w} and, for $\vert t\vert > C$,
$\Gamma(w) \ll \vert t \vert^{-1} \exp(-\frac{\pi}{2}\vert t \vert)$
by \eqref{Stirling}.
Letting $L_{C}=\{-1/2+it: \vert t \vert > C\}$  we have
\[
\int_{L_{C}} 
\frac{\zeta'}{\zeta}(\ell w)  \Gamma(w) z^{-w} \, \dx w
  \ll_\ell
   \vert z \vert^{1/2}
   \int_{L_{C}} 
     \frac{\log \vert t \vert}{\vert t \vert}
     \exp
     \Bigl(
     -\frac{\pi}{2}\vert t \vert + t \arctan(  y / a )
     \Bigr)
     \, \dx t.
\]
If $t y\leq 0$ we call $\eta$ the quantity 
$\frac{\pi}{2}+ \vert\arctan(  y / a )\vert \in [\pi/2, \pi)$.
If $\vert y \vert \leq a$ we define $\eta$ as 
$\frac{\pi}{2} - \arctan( y / a )>\frac{\pi}{2} - \arctan(1)=
\frac{\pi}{4}$.
In the remaining case ($\vert y \vert > a$ and $ty > 0$)
we set $\eta=\arctan(a/ \vert y\vert) \gg a/\vert y \vert$.
Now fix $C$ such that $C\eta<1$ (e.g., $C=1/\pi$ is allowed). 
Letting $u = \eta t$, we get
\begin{align}
  \notag
  \int_{L_{C}}
  \frac{\zeta'}{\zeta}(\ell w)  \Gamma(w) z^{-w} \, \dx w
  &\ll_\ell
  \vert z \vert^{1/2}
  \int_C^{+\infty} e^{-\eta t} \, \frac{\log t}t \, \dx t
  =
  \vert z \vert^{1/2}
  \int_{C \eta}^{+\infty} e^{-u} \, \frac{\log(u / \eta)}u \, \dx u
  \\
  \notag
  &=
  \vert z \vert^{1/2}
  \int_{C \eta}^{+\infty} e^{-u} \, \frac{\log u}u \, \dx u
  +
  \vert z \vert^{1/2}
  \log(1 / \eta)
  \int_{C \eta}^{+\infty} e^{-u} \frac{\dx u}u \\
\label{splitting}
  &=
  J_1 + J_2.
\end{align}
We remark that $0 \le u^{-1} \log u \le e^{-1}$ for $u \ge 1$, since
the maximum of $u^{-1} \log u$ is attained at $u = e$.
Since
\[
  0
  \le
  \int_1^{+\infty} e^{-u} \, \frac{\log u}u \, \dx u
  \le
  e^{-1}
  \int_1^{+\infty} e^{-u} \, \dx u
  \ll
  1
\]
and
\[
  \Bigl\vert
    \int_{C \eta}^1 e^{-u} \, \frac{\log u}u \, \dx u
  \Bigr\vert
  \le
  \int_{C \eta}^1 \frac{-\log u}u \, \dx u
  =
  \Bigl[ - \frac12 \log^2 u
  \Bigr]_{C \eta}^1
  \ll
  \log^2 (1 / \eta)
\]
we have that $J_{1} \ll \vert z \vert^{1/2} \log^2 (1 / \eta)$.
For $J_{2}$ it is sufficient to remark that
\[
  0
  \le
  J_2
  \le
  \vert z \vert^{1/2}
  \log(1 / \eta)
  \Bigl(
    \int_{C \eta}^1 \frac{\dx u}u
    +
    \int_1^{+\infty} e^{-u} \dx u
  \Bigr)
  \ll
  \vert z \vert^{1/2}
  \log^2 (1 / \eta).
\]
Inserting the last two estimates in \eqref{splitting}, recalling the
definition of $\eta$, remarking that the integration over
$\vert t \vert \leq C$ gives immediately a contribution $\ll_\ell \vert z \vert^{1/2}$,
we get that 
\[
  \int_{\mathcal{L}_\ell} 
  \frac{\zeta'}{\zeta}(\ell w) \Gamma(w)z^{-w} \, \dx w
\ll_\ell 
  \vert z \vert^{1/2}
  \begin{cases}
    1 & \text{if $\vert y \vert \leq a$} \\
    1 +\log^2 (\vert y\vert/a) & \text{if $\vert y \vert > a$.}
  \end{cases}
\]
provided that $4\nmid \ell$.
Recalling $(\zeta'/\zeta) (0) = \log (2\pi)$ and remarking that
\begin{equation}
\label{sum-trivial-zeros-estim}
  \sum_{m=1}^{\ell/4} \Gamma\Bigl(-\frac{2m}{\ell}\Bigr) z^{2m/\ell} \ll_\ell  \vert z \vert^{1/2},
\end{equation}
we obtain that the case $4\nmid \ell$ of the  lemma is proved.

Assume now that $4\mid \ell$. The computation over $L_C$
can be performed as in the previous case; we can also choose
$C=1/\pi$ as we did before.
On the vertical segments $\mathcal{S}$ given by $\Re(w)=-1/2$, $\vert \Im(w) \vert \in[1/\ell, C]$, we exploit  
the boundedness of the $\Gamma$-function and the estimate $\vert z^{-w}\vert \ll \vert z \vert^{1/2}$
which holds on $\mathcal{S}$ since the argument of $z$ is bounded there. 
This gives 
\[
  \frac{1}{2\pi i}
  \int_{\mathcal{S}} 
  \frac{\zeta'}{\zeta}(\ell w) \Gamma(w)z^{-w} \, \dx w
  \ll_\ell 
  \vert z \vert^{1/2}.
\]
It remains to consider the contribution over 
$\gamma_{\ell}$; on this path we can again make use 
of the boundedness of the $\Gamma$-function and 
 that  $\vert z^{-w}\vert \ll \vert z \vert^{1/2}$
since the argument of $z$ is bounded on $\gamma_{\ell}$. 
This leads to 
\[
  \frac{1}{2\pi i}
  \int_{\gamma_{\ell}} 
  \frac{\zeta'}{\zeta}(\ell w) \Gamma(w)z^{-w} \, \dx w
  \ll_\ell 
  \vert z \vert^{1/2}.
\]
Summing up, for $4 \mid \ell$  we obtain that the integral in \eqref{Mellin2} is dominated by the right
hand side of \eqref{expl-form-err-term-strong}
and this, together with  \eqref{sum-trivial-zeros-estim} and   $ (\zeta'/\zeta) (0) = \log (2\pi)$, 
proves this case of the  lemma.
\end{Proof}

We remark that, at the cost of some other complications in the details, 
Lemma \ref{Linnik-lemma-gen} can be extended to the case $\ell\in \R$, $\ell>0$.

In the next sections we will need to perform several times a set of 
similar computations; we collected them in the following two lemmas,
which extend Lemmas~4.2 and~4.3 in \cite{LanguascoZ2015a}.

\begin{Lemma}
\label{series-int-zeros}
Let $\ell\ge 1$ be an integer, let $\beta + i \gamma$ run over the non-trivial zeros of the Riemann
zeta-function and $\alpha > 1$ be a parameter.
For any fixed $c \ge 0$ the series
\[
  \sum_{\rho \colon \gamma > 0}
  \Bigr(\frac{\gamma}{\ell}\Bigl)^{\beta/\ell-1/2}
    \int_1^{+\infty} (\log u)^c
      \exp\Bigl( - \frac{\gamma}{\ell} \arctan\frac 1u \Bigr)
      \frac{\dx u}{u^{\alpha+\beta/\ell}}
\]
converges provided that $\alpha > 3/2$.
For $\alpha \le 3/2$ the series does not converge.
\end{Lemma}

\begin{Proof} 
Setting $y = \arctan(1 / u)$, for any real $\gamma > 0$ we have
\begin{align*}
  \int_1^{+\infty} \exp\Bigl( - \frac{\gamma}{\ell} \arctan\frac 1u \Bigr)
    \frac{\dx u}{u^{\alpha+\beta/\ell}} 
  &=
  \int_0^{\pi / 4}
    \exp\Bigl(-\frac{\gamma y}{\ell}   \Bigr) \,
    \frac{(\sin y)^{\alpha+\beta/\ell - 2}}{(\cos y)^{\alpha+\beta/\ell}} \, \dx y\\
  &\ll_\alpha
  \int_0^{\pi / 4}
    \exp\Bigl(-\frac{\gamma y}{\ell}\Bigr ) \, y^{\alpha+\beta/\ell - 2} \, \dx y \\
  &=
  \Bigl( \frac{\gamma}{\ell} \Bigr)^{1 - \alpha - \beta/\ell}
  \int_0^{\pi  \gamma / (4\ell)}
    \exp(-w) \, w^{\alpha+\beta/\ell - 2} \, \dx w \\
  &\ll_{\alpha,\ell}
  \Bigr(\frac{\gamma}{\ell}\Bigl)^{1 - \alpha - \beta/\ell} \,
  \bigl(\Gamma(\alpha-1)+\Gamma(\alpha+1/\ell -1)\bigr),
\end{align*}
since $0 < \beta < 1$.
This shows that the series over $\gamma$ converges for $\alpha > 3/2$.
For $\alpha = 3/2$ essentially the same computation shows that the
integral is $\gg \gamma^{-1/2 - \beta/\ell}$ and it is well known that in
this case the series over zeros diverges.
\end{Proof}

\begin{Lemma}
\label{series-int-zeros-alt-sign}
Let $\ell\ge 1$ be an integer,   $\alpha > 1$, $z=a+iy$, $a\in(0,1)$ and $y\in \R$.
Let further $\rho=\beta+i\gamma$ run over the non-trivial zeros of
the Riemann zeta-function.  We have
\[
  \sum_{\rho}
    \Bigl\vert \frac{\gamma}{\ell}\Bigr\vert ^{\beta/\ell-1/2}
    \int_{\Y_1 \cup \Y_2} \exp\Bigl(\frac{\gamma }{\ell} \arctan\frac{y}{a} - \frac\pi2\Bigl \vert \frac{\gamma}{\ell}  \Bigr\vert\Bigr)
      \frac{\dx y}{\vert z \vert ^{\alpha+\beta/\ell}}
  \ll_{\alpha,\ell}
  a^{1-\alpha-1/\ell},
\]
where $\Y_1=\{y\in \R\colon y\gamma \leq 0\}$ and 
$\Y_2=\{y\in [-a,a] \colon y\gamma > 0\}$.
The result remains true if we insert in the integral a factor
$(\log (\vert y\vert /a))^c$, for any fixed $c \ge 0$.
\end{Lemma}

\begin{Proof}
We first work on $\Y_1$.
By symmetry, we may assume that $\gamma > 0$.
For $y \in(-\infty, 0]$ we have
$(\gamma/\ell) \arctan(y/a) -\frac \pi2 \vert \gamma/\ell\vert  \le - \frac \pi2 \vert \gamma/\ell\vert $
and hence the quantity we are estimating becomes
\begin{align*}
  \sum_{\rho \colon \gamma > 0} 
  \Bigl( \frac{\gamma}{\ell} \Bigr)^{\beta/\ell-1/2}
  \exp\Bigl( -\frac \pi2 \frac{\gamma}{\ell} \Bigr)
    \int_{-\infty}^0 \frac{\dx y}{\vert z \vert ^{\alpha+\beta/\ell}}
 & 
 \ll_{\alpha,\ell}
  \sum_{\rho \colon \gamma > 0} 
  \Bigr(\frac{\gamma}{\ell}\Bigl)^{\beta/\ell-1/2}
  \exp\Bigl( -\frac \pi2 \frac{\gamma}{\ell}  \Bigr)
  a^{1-\alpha-\beta/\ell}
 \\&
  \ll_{\alpha,\ell}
  a^{1-\alpha-1/\ell},
\end{align*}
using $0<\beta<1$, standard zero-density estimates and \eqref{z^-1}.
We consider now the integral over $\Y_2$.
Again by symmetry we can assume that $\gamma > 0$ and so we get
\begin{align*}
  \sum_{\rho \colon \gamma > 0}
   \Bigl( \frac{\gamma}{\ell} \Bigr)^{\beta/\ell-1/2} &
    \int_0^{a} \exp\Bigl( \frac{\gamma}{\ell} \bigl( \arctan \frac{y}{a} -\frac \pi2 \bigr) \Bigr)
     \frac{\dx y}{\vert z \vert ^{\alpha+\beta/\ell}}
  \\& \ll
  \sum_{\rho \colon \gamma > 0}
   \Bigl( \frac{\gamma}{\ell} \Bigr)^{\beta/\ell-1/2} 
  \exp\Bigl( -\frac \pi4 \frac{\gamma}{\ell}  \Bigr)
  \int_0^{a} \frac{\dx y}{\vert z \vert ^{\alpha+\beta/\ell}} \\
  &  \ll_{\alpha,\ell}
  \sum_{\rho \colon \gamma > 0}
   \Bigr(\frac{\gamma}{\ell}\Bigl)^{\beta/\ell-1/2} \exp\Bigl(  -\frac \pi4 \frac{\gamma}{\ell}  \Bigr)
  a^{1-\alpha-\beta/\ell} 
    \ll_{\alpha,\ell}
  a^{1-\alpha-1/\ell}
\end{align*}
arguing as above.
The other assertions are proved in the same way.
\end{Proof}

\section{Interchange of the series over zeros with the line integral in $I_3$}
\label{first-exchange}

We need $k>1/2$ in this section. We need to establish the convergence of
\begin{equation}
\label{conv-integral}
  \sum_{\rho}
    \Bigl\vert \Gamma\Bigr(\frac{\rho}{\ell}\Bigl) \Bigr\vert
    \int_{(\frac{1}{N})} \vert e^{N z} \vert \, \vert z \vert ^{- k - 1} \,
      \vert z^{-\rho/\ell} \vert \, \vert \dx z \vert.
\end{equation}
By \eqref{z^w} and the Stirling formula \eqref{Stirling}, we are left
with estimating
\begin{equation}
\label{conv-integral-1}
  \sum_{\rho}
    \Bigl\vert \frac{\gamma}{\ell}\Bigr \vert ^{\beta/\ell_j - 1/2}
    \int_{\R} \exp\Bigl( \frac{\gamma}{\ell} \arctan(N y) -\frac \pi2 \Bigl \vert \frac{\gamma}{\ell} \Bigr\vert \Bigr)
             \frac{\dx y}{\vert z \vert ^{k + 1 +\beta/\ell}}. 
\end{equation}
We have just to consider the case $\gamma y >0$, $\vert y \vert > 1/N$
since in the other cases the total contribution is  $\ll_{k,\ell}N^{k + 1/\ell}$
by Lemma \ref{series-int-zeros-alt-sign} with $\alpha=k+1$ and $a=1/N$.
By symmetry, we may assume that $\gamma > 0$.
We have that the integral in \eqref{conv-integral-1} is
\begin{align*}  
  &\ll_{\ell}
  \sum_{\rho \colon \gamma > 0}
    \Bigr(\frac{\gamma}{\ell}\Bigl) ^{\beta/\ell - 1/2}
    \int_{1 / N}^{+\infty} \exp\Bigl( -  \frac{\gamma}{\ell} \arctan\frac 1{N y} \Bigr)
      \frac{\dx y}{y^{k + 1 +\beta/\ell}} \\
  &=
  N^{k}
  \sum_{\rho \colon \gamma > 0}
    N^{\beta/\ell}
    \Bigr(\frac{\gamma}{\ell}\Bigl) ^{\beta/\ell - 1/2}
    \int_1^{+\infty} \exp\Bigl( - \frac{\gamma}{\ell} \arctan\frac 1u \Bigr)
      \frac{\dx u}{u^{k + 1 +\beta/\ell}}.
\end{align*}
For $k > 1 / 2$ this is $\ll_{k,\ell} N^{k + 1/\ell}$ by
Lemma~\ref{series-int-zeros}. This implies
that the integrals in \eqref{conv-integral-1} and in \eqref{conv-integral} 
are both  $\ll_{k,\ell} N^{k + 1/\ell}$ and
hence the exchange steps for $I_3$ are fully justified.

\section[Interchange of the series over zeros]
{Interchange of the series over zeros with the line integral in $I_4$}
\label{first-exchange-bis}

We need $k>1/2-1/\ell_2$ in this section. We need to establish the convergence of
\begin{equation}
\label{conv-integral-bis}
  \sum_{\rho}
    \Bigl\vert \Gamma\Bigr(\frac{\rho}{\ell_{1}}\Bigl) \Bigr\vert
    \int_{(\frac{1}{N})} \vert e^{N z} \vert \, \vert z \vert ^{- k - 1-1/\ell_2} \,
      \vert z^{-\rho/\ell_1} \vert \, \vert \dx z \vert
\end{equation}
and of the case in which $\ell_1$ and $\ell_2$ are interchanged.
By \eqref{z^w} and the Stirling formula \eqref{Stirling}, we are left
with estimating
\begin{equation}
\label{conv-integral-bis-1}
  \sum_{\rho}
    \Bigl\vert \frac{\gamma}{\ell_1}\Bigr\vert ^{\beta/\ell_1 - 1/2}
    \int_{\R} \exp\Bigl( \frac{\gamma}{\ell_1} \arctan(N y) -\frac \pi2 \Bigl \vert \frac{\gamma}{\ell_1} \Bigr\vert\Bigr)
             \frac{\dx y}{\vert z \vert ^{k + 1 + 1/\ell_2+\beta/\ell_1}}. 
\end{equation}
We have just to consider the case $\gamma y >0$, $\vert y \vert > 1/N$
since in the other cases the total contribution is  $\ll_{k,\ell_1,\ell_2} N^{k + \lambda}$
by Lemma \ref{series-int-zeros-alt-sign} with $\alpha=k+1+ 1/\ell_2$ and $a=1/N$.
By symmetry, we may assume that $\gamma > 0$.
We have that the integral in \eqref{conv-integral-bis-1} is
\begin{align*}  
  &\ll_{\ell_1}
  \sum_{\rho \colon \gamma > 0}
    \Bigr(\frac{\gamma}{\ell_1}\Bigl) ^{\beta/\ell_1- 1/2}
    \int_{1 / N}^{+\infty} \exp\Bigl( -  \frac{\gamma}{\ell_1} \arctan\frac 1{N y} \Bigr)
      \frac{\dx y}{y^{k + 1 + 1/\ell_2 +\beta/\ell_1}} \\
  &=
  N^{k+1/\ell_2}
  \sum_{\rho \colon \gamma > 0}
    N^{\beta/\ell_1}
    \Bigr(\frac{\gamma}{\ell_1}\Bigl) ^{\beta/\ell_1- 1/2}
    \int_1^{+\infty} \exp\Bigl( - \frac{\gamma}{\ell_1} \arctan\frac 1u \Bigr)
      \frac{\dx u}{u^{k + 1 + 1/\ell_2+\beta/\ell_1}}.
\end{align*}
For $k > 1/2-1/\ell_2$ this is $\ll_{k,\ell_1,\ell_2} N^{k + \lambda}$ by
Lemma~\ref{series-int-zeros}. This implies
that the integrals in \eqref{conv-integral-bis-1} and in \eqref{conv-integral-bis} 
are both  $\ll_{k,\ell_1,\ell_2} N^{k + \lambda}$ and
hence the exchange step for $I_4$  is fully justified.

\section[Interchange of the double series over zeros]
{Interchange of the double series over zeros with the line integral in $I_5$}
\label{exchange-double-sum-rhos}

 We need $k>1$ in this section.
Arguing as in Sections  \ref{first-exchange}-\ref{first-exchange-bis},
we first need to establish the
convergence of
\begin{equation}
\label{conv-integral-2}
  \sum_{\rho_{1}}
    \Bigl\vert \Gamma\Bigr(\frac{\rho_1}{\ell_{1}}\Bigl) \Bigr\vert  
    \int_{(\frac{1}{N})}  
      \Bigl\vert \sum_{\rho_{2}} \Gamma\Bigr(\frac{\rho_2}{\ell_{2}}\Bigl) z^{- \rho_{2}/\ell_2} \Bigr\vert
      \vert e^{N z} \vert \,
      \vert z \vert ^{- k - 1} \, \vert z^{- \rho_{1}/\ell_1} \vert \, \vert \dx z \vert.
\end{equation}

Using the Prime Number Theorem and 
\eqref{expl-form-err-term-strong},
 we first remark
that
\begin{equation} 
\label{sum-over-rho-new}
 \Bigl\vert \sum_{\rho_{2}} \Gamma\Bigr(\frac{\rho_2}{\ell_{2}}\Bigl) z^{- \rho_{2}/\ell_2} \Bigr\vert
    \ll_{\ell_2}
  N^{1/\ell_2} + \vert z\vert ^{1/2}
    \log^2 (2N\vert y\vert). 
  \end{equation}
By symmetry, we may assume that $\gamma_{1}> 0$.
By \eqref{sum-over-rho-new}, 
\eqref{z^-1}, \eqref{z^w} and \eqref{lambda-def}, 
for $y \in(-\infty, 0]$  we are first led to estimate
\begin{align*}
  \sum_{\rho_{1} \colon \gamma_{1} > 0} 
    \Bigl(\frac{\gamma_{1}}{\ell_1}\Bigr)^{\beta_{1}/\ell_1 - 1/2}
    \exp\Bigl( -\frac \pi2  \frac{\gamma_{1}}{\ell_1}  \Bigr)
  &\Bigl(
    \int_{-1/N}^{0}  N^{ k +  1+1/\ell_2 + \beta_{1}/\ell_1} \, \dx y
    +
    N^{1/\ell_2}\int_{-\infty}^{-1/N}  
     \frac{\dx y}{\vert y\vert^{k + 1+\beta_{1}/\ell_1}}    
    \\
  &+
  \int_{-\infty}^{-1/N}  \log^{2} (2N \vert y \vert)
    \frac{\dx y}{\vert y\vert^{k + 1/2+\beta_{1}/\ell_1}}     
  \Bigr)
  \ll_{k,\ell_1,\ell_2}
  N^{k + \lambda}
\end{align*}
by the same argument used in the proof of Lemma
\ref{series-int-zeros-alt-sign} 
with $\alpha=k+1/2$ and $a=1/N$.
On the other hand, for $y > 0$ we split the range of integration into
$(0, 1 / N] \cup (1 / N, +\infty)$.
By \eqref{sum-over-rho-new}, \eqref{z^-1} and Lemma \ref{series-int-zeros-alt-sign}
with $\alpha=k+1$ and $a=1/N$, on  $[0,1/N]$ we have
\begin{align*}
 N^{1/\ell_2} \sum_{\rho_{1} \colon \gamma_{1} > 0}
    \Bigl(\frac{\gamma_{1}}{\ell_1}\Bigr)^{\beta_{1}/\ell_1 - 1/2}
    \int_0^{1 / N} \exp\Bigl(  \frac{\gamma_{1}}{\ell_1}\bigl( \arctan(N y) -\frac \pi2 \bigr) \Bigr)
    \frac{\dx y}{\vert z\vert^{k + 1+\beta_{1}/\ell_1}}   
  &\ll_{k,\ell_1,\ell_2}
   N^{k + \lambda}.
\end{align*}
On the other interval, again by \eqref{z^-1}, 
we have to estimate
\begin{align*}
  &
  \sum_{\rho_{1} \colon \gamma_{1} > 0}
   \Bigl(\frac{\gamma_{1}}{\ell_1}\Bigr)^{\beta_{1}/\ell_1 - 1/2}
    \int_{1 / N}^{+\infty} \exp\Bigl( - \frac{\gamma_{1}}{\ell_1}\arctan\frac 1{N y} \Bigr) 
      \frac{ N^{1/\ell_2} + y^{1/2}\log^{2} (2Ny)}{y^{k + 1+\beta_{1}/\ell_1}} \dx y
 \\
  &=
  N^{k}
  \sum_{\rho_{1} \colon \gamma_{1} > 0} \!\!\!\!
  N^{\beta_{1}/\ell_1 }
   \Bigl(\frac{\gamma_{1}}{\ell_1}\Bigr)^{\beta_{1}/\ell_1 - 1/2}
    \int_{1 }^{+\infty} \exp\Bigl( -\frac{\gamma_{1}}{\ell_1}\arctan\frac 1{u} \Bigr)
      \frac{N^{1/\ell_2}  + u^{1/2}N^{-1/2}\log^{2} (2u)}{u^{k + 1+\beta_{1}/\ell_1}} \dx u.
\end{align*}
Recalling \eqref{lambda-def},
Lemma~\ref{series-int-zeros} 
with $\alpha=k+1/2$ shows that the last
term is $\ \ll_{k,\ell_1,\ell_2} N^{k + \lambda}$.
This implies that the integral in 
\eqref{conv-integral-2} 
is $\ll_{k,\ell_1,\ell_2} N^{k + \lambda}$ 
provided that $k > 1$ and hence we can exchange the first summation 
with the integral in this case.

To exchange the second summation we have to consider
\begin{equation}
\label{conv-integral-3}
  \sum_{\rho_{1}}
    \Bigl\vert \Gamma\Bigr(\frac{\rho_1}{\ell_{1}}\Bigl)\Bigr\vert  
  \sum_{\rho_{2}}
 \Bigl\vert \Gamma\Bigr(\frac{\rho_2}{\ell_{2}}\Bigl)\Bigr\vert    
      \int_{(\frac{1}{N})} \vert 
    e^{N z}\vert  \vert z\vert^{- k - 1} 
    \vert z^{- \rho_{1}/\ell_1} \vert 
    \vert z^{- \rho_{2}/\ell_1} \vert \, \vert \dx z \vert.
\end{equation}
By symmetry, we can consider  $\gamma_{1},\gamma_{2}> 0$ or
$\gamma_{1} >0$, $\gamma_{2}< 0$.

Assuming $\gamma_{1},\gamma_{2}> 0$, 
for $y \leq  0$ we have
$(\gamma_{j}/\ell_j) \arctan(N y) -\frac \pi2  (\gamma_{j}/\ell_j) \le - \frac \pi2 (\gamma_{j}/\ell_j)$,
$j=1,2$, and, by \eqref{z^w}, the corresponding contribution to
\eqref{conv-integral-3} is $ \ll_{k,\ell_1,\ell_2} N^{k + \lambda}$
since
\begin{align*}
  \sum_{\rho_{1} \colon \gamma_{1} > 0}
   & \Bigl(\frac{\gamma_{1}}{\ell_1}\Bigr)^{\beta_{1}/\ell_1 - 1/2}
    \exp\Bigl( -\frac \pi2  \frac{\gamma_{1}}{\ell_1}  \Bigr)
  \sum_{\rho_{2} \colon \gamma_{2} > 0} 
   \Bigl(\frac{\gamma_{2}}{\ell_2}\Bigr)^{\beta_{2}/\ell_2 - 1/2}\exp\Bigl( -\frac \pi2  \frac{\gamma_{2}}{\ell_2}  \Bigr)
  \Bigl(
    \int_{-\infty}^{0} 
     \frac{\dx y}{\vert z\vert^{k + 1+\beta_{1}/\ell_1+\beta_{2}/\ell_2}}   
    \Bigr)
    \\
  &\ll_k
  N^{k + \lambda}
   \sum_{\rho_{1} \colon \gamma_{1} > 0} 
   \Bigl(\frac{\gamma_{1}}{\ell_1}\Bigr)^{\beta_{1}/\ell_1 - 1/2} \exp\Bigl( -\frac \pi2  \frac{\gamma_{1}}{\ell_1}  \Bigr)
   \sum_{\rho_{2} \colon \gamma_{2} > 0} 
       \Bigl(\frac{\gamma_{2}}{\ell_2}\Bigr)^{\beta_{2}/\ell_2 - 1/2}\exp\Bigl( -\frac \pi2  \frac{\gamma_{2}}{\ell_2}  \Bigr),
\end{align*}
using standard zero-density estimates, \eqref{z^-1} and \eqref{lambda-def}.
On the other hand, for $y > 0$ we split the range of integration into
$(0, 1 / N] \cup (1 / N, +\infty)$.
On the first interval we have
\begin{align*}
&
  \sum_{\rho_{1} \colon \gamma_{1} > 0} 
   \Bigl(\frac{\gamma_{1}}{\ell_1}\Bigr)^{\beta_{1}/\ell_1 - 1/2} 
  \sum_{\rho_{2} \colon \gamma_{2} > 0}
 \Bigl(\frac{\gamma_{2}}{\ell_2}\Bigr)^{\beta_{2}/\ell_2 - 1/2} 
 \\& 
 \hskip1cm
 \times
  \int_0^{1 / N}
    \exp\Bigl( \Bigl(\frac{\gamma_{1}}{\ell_1}+\frac{\gamma_{2}}{\ell_2}\Bigr) \bigl( \arctan(N y) -\frac \pi2 \bigr) \Bigr)
       \frac{\dx y}{\vert z\vert^{k + 1+\beta_{1}/\ell_1+\beta_{2}/\ell_2}}    \\
  &\ll
  \sum_{\rho_{1} \colon \gamma_{1} > 0}
  \Bigl(\frac{\gamma_{1}}{\ell_1}\Bigr)^{\beta_{1}/\ell_1 - 1/2} 
  \sum_{\rho_{2} \colon \gamma_{2} > 0} 
    \Bigl(\frac{\gamma_{2}}{\ell_2}\Bigr)^{\beta_{2}/\ell_2 - 1/2} 
   \exp\Bigl( -\frac \pi4  \Bigl(\frac{\gamma_{1}}{\ell_1}+\frac{\gamma_{2}}{\ell_2}\Bigr)  \Bigr)
  \int_0^{1 / N}  N^{k + 1 +\beta_{1}/\ell_1 +\beta_{2}/\ell_2} \, \dx y \\
  &\ll_{k,\ell_1,\ell_2}
    N^{k + \lambda}
  \sum_{\rho_{1} \colon \gamma_{1} > 0} 
  \Bigl(\frac{\gamma_{1}}{\ell_1}\Bigr)^{\beta_{1}/\ell_1 - 1/2}  \exp\Bigl(  -\frac \pi4 \frac{\gamma_{1}}{\ell_1}  \Bigr)
  \sum_{\rho_{2} \colon \gamma_{2} > 0} 
  \Bigl(\frac{\gamma_{2}}{\ell_2}\Bigr)^{\beta_{2}/\ell_2 - 1/2}  \exp\Bigl( -\frac \pi4 \frac{\gamma_{2}}{\ell_2}  \Bigr),
\end{align*}
which is also $\ll_{k,\ell_1,\ell_2} N^{k + \lambda}$, arguing as
above. With similar computations, on the other interval we have
\begin{align*}
  \sum_{\rho_{1} \colon \gamma_{1} > 0} &
  \Bigl(\frac{\gamma_{1}}{\ell_1}\Bigr)^{\beta_{1}/\ell_1 - 1/2} 
  \sum_{\rho_{2} \colon \gamma_{2} > 0}
  \Bigl(\frac{\gamma_{2}}{\ell_2}\Bigr)^{\beta_{2}/\ell_2 - 1/2} 
  \\& \hskip3cm
 \times
  \int_{1 / N}^{+\infty}
    \exp\Bigl( \Bigl(\frac{\gamma_{1}}{\ell_1}+\frac{\gamma_{2}}{\ell_2}\Bigr) \bigl( \arctan(N y) -\frac \pi2 \bigr) \Bigr)
    \frac{\dx y}{y^{k + 1+\beta_{1}/\ell_1+\beta_{2}/\ell_2}} \\
  &=
  N^{k}
  \sum_{\rho_{1} \colon \gamma_{1} > 0}
  N^{\beta_{1}/\ell_1}
 \Bigl(\frac{\gamma_{1}}{\ell_1}\Bigr)^{\beta_{1}/\ell_1 - 1/2} 
  \sum_{\rho_{2} \colon \gamma_{2} > 0}
  N^{\beta_{2}/\ell_2}
 \Bigl(\frac{\gamma_{2}}{\ell_2}\Bigr)^{\beta_{2}/\ell_2 - 1/2} 
 \\& \hskip3cm
 \times
    \int_{1}^{+\infty} \exp\Bigl(- \Bigl(\frac{\gamma_{1}}{\ell_1}+\frac{\gamma_{2}}{\ell_2}\Bigr) \arctan\frac 1u \Bigr)
      \frac{\dx u}{u^{k + 1+\beta_{1}/\ell_1+\beta_{2}/\ell_2}}.      
      \end{align*}
Arguing as in the proof of Lemma~\ref{series-int-zeros}, we prove that the integral
on the right is $\sameorder_{k,\ell_1,\ell_2} (\gamma_1 + \gamma_2)^{-k-\beta_1/\ell_1-\beta_2/\ell_2}$.
The inequality
\begin{equation}
\label{ineq-zeros}
  \frac{\gamma_1^{\beta_1/\ell_1 - 1/2} \gamma_2^{\beta_2/\ell_2- 1/2}}
       {(\gamma_1 + \gamma_2)^{\beta_1/\ell_1+\beta_2/\ell_2}}
  \le
  \frac1{\gamma_1^{1/2} \gamma_2^{1/2}}
\end{equation}
shows, using \eqref{lambda-def}, that it is sufficient to consider
\begin{align*}
  N^{k} &
  \sum_{\rho_{1} \colon \gamma_{1} > 0}
  \sum_{\rho_{2} \colon \gamma_{2} > 0}
    N^{\frac{\beta_{1}}{\ell_1}+\frac{\beta_{2}}{\ell_2}}
    \frac1{\gamma_{1} ^{1/2}\gamma_{2} ^{1/2}(\gamma_1 + \gamma_2)^{k}} \\
 &\ll_{k,\ell_1,\ell_2}
  N^{k + \lambda}
  \sum_{\rho_{1} \colon \gamma_{1} > 0}
    \frac1{\gamma_1^{k + 1/2}}
    \sum_{\rho_{2} \colon  0< \gamma_2 \le \gamma_1}
      \frac1{\gamma_2^{1/2}}
  \ll_{k,\ell_1,\ell_2}
  N^{k + \lambda}
  \sum_{\rho_{1} \colon \gamma_{1} > 0}
   \frac{\log \gamma_1}{\gamma_1^{k}}
\end{align*}
and the last series over zeros converges for $k > 1$.
 
Assume now $\gamma_{1}>0$, $\gamma_{2}< 0$.
For $y \leq  0$ we have
$\frac{\gamma_{1}}{\ell_1}\arctan(N y) -\frac \pi2  \frac{\gamma_{1}}{\ell_1}  
\le - \frac \pi2  \frac{\gamma_{1}}{\ell_1}$,  
by \eqref{z^-1} and \eqref{lambda-def} the corresponding contribution to \eqref{conv-integral-3} is
\begin{align*} 
&\ll_{k,\ell_1,\ell_2} 
  \sum_{\rho_{1} \colon \gamma_{1} > 0}  
   \Bigl(\frac{\gamma_{1}}{\ell_1}\Bigr)^{\beta_{1}/\ell_1 - 1/2} 
   \exp\Bigl( -\frac \pi2  \frac{\gamma_{1}}{\ell_1}  \Bigr)
   \\& \times
  \Bigl\{
   \sum_{\rho_{2} \colon \gamma_{2} < 0}   
     \Bigl\vert \frac{\gamma_{2}}{\ell_2}  \Bigr\vert ^{\beta_{2}/\ell_2 - 1/2} 
     \Bigl[
     \exp\Bigl( -\frac \pi4 \Bigl\vert \frac{\gamma_{2}}{\ell_2}  \Bigr\vert   \Bigr)
    \int_{-1/N}^{0}
    N^{k + 1 +\beta_{1}/\ell_1+\beta_{2}/\ell_2} \, \dx y 
    \\
    &
    \hskip 1cm
    +    
    \int_{-\infty}^{-1/N} 
    \exp\Bigl(- \Bigl\vert \frac{\gamma_{2}}{\ell_2}  \Bigr\vert \bigl( \arctan(N y) +\frac \pi2 \bigr) \Bigr)
    \frac{\dx y}{\vert y\vert ^{k + 1+\beta_{1}/\ell_1+\beta_{2}/\ell_2}}  
     \Bigr]
  \Bigr\}  \\
  &\ll_{k,\ell_1,\ell_2} 
  N^{k + \lambda}
   \sum_{\rho_{1} \colon \gamma_{1} > 0}  
   \Bigl(\frac{\gamma_{1}}{\ell_1}\Bigr)^{\beta_{1}/\ell_1 - 1/2} 
   \exp\Bigl( -\frac \pi2  \frac{\gamma_{1}}{\ell_1}  \Bigr)
\sum_{\rho_{2} \colon \gamma_{2} < 0} 
\Bigl\vert \frac{\gamma_{2}}{\ell_2}  \Bigr\vert^{\beta_{2}/\ell_2 - 1/2}  
 \exp\Bigl( -\frac \pi4\Bigl\vert \frac{\gamma_{2}}{\ell_2}  \Bigr\vert \Bigr)
\\
&
\hskip 1cm +  
  N^{k + \lambda}
 \!\!    \sum_{\rho_{1} \colon \gamma_{1} > 0}  
   \Bigl(\frac{\gamma_{1}}{\ell_1}\Bigr)^{\beta_{1}/\ell_1 - 1/2} 
   \exp\Bigl( -\frac \pi2  \frac{\gamma_{1}}{\ell_1}  \Bigr)
 \!\! \sum_{\rho_{2} \colon \gamma_{2} < 0} \!\!
 \Bigl\vert \frac{\gamma_{2}}{\ell_2}  \Bigr\vert^{\beta_{2}/\ell_2 - 1/2} 
  \\& 
  \hskip 2cm\times
  \int_1^{+\infty} \!\!\!\!
   \exp\Bigl( -\Bigl\vert \frac{\gamma_{2}}{\ell_2}  \Bigr\vert\arctan\frac 1u \Bigr)
      \frac{\dx u}{u^{k + 1 +\beta_{1}/\ell_1+\beta_{2}/\ell_2}} 
    \\
  &\ll_{k,\ell_1,\ell_2} 
  N^{k + \lambda}
  +
 N^{k + \lambda}
  \sum_{\rho_{1} \colon \gamma_{1} > 0} 
  \Bigl(\frac{\gamma_{1}}{\ell_1}\Bigr)^{\beta_{1}/\ell_1 - 1/2} 
  \exp\Bigl( -\frac \pi2 \gamma_{1} \Bigr) 
 \ll_{k,\ell_1,\ell_2} 
  N^{k + \lambda}
\end{align*}
for $k > 1/2$, by Lemma~\ref{series-int-zeros} and standard zero-density
estimates.

On the other hand, the case $\gamma_{1}>0$, $\gamma_{2}< 0$ and $y > 0$
can be estimated in a similar way essentially exchanging the role of
$\gamma_{1}$ and $\gamma_{2}$ in the previous argument.

This implies that the integral in \eqref{conv-integral-3} is $ \ll_{k,\ell_1,\ell_2} N^{k + \lambda}$ 
provided that $k > 1$. Combining the convergence conditions for 
\eqref{conv-integral-2}-\eqref{conv-integral-3},
we see that we can exchange both summations with the integral provided that $k>1$.
 
\section{Convergence of the double sum over zeros}
\label{sec:double-sum}

In this section we prove that the double sum on the right of
\eqref{def-M5} converges absolutely 
for every $k > 1 / 2$; the other series in \eqref{def-M3} and
\eqref{def-M4} clearly converge for $k>0$ or better.
We need \eqref{Stirling} uniformly for $x \in [0, k + 3]$ and
$\vert y \vert \ge T$, where $T$ is large but fixed: this provides
both an upper and a lower bound for $\vert \Gamma(x + i y) \vert$.
Let
\begin{equation*} 
  \Sigma 
  =
  \sum_{\rho_1} \sum_{\rho_2}
    \Bigl\vert \frac{\Gamma(\rho_1/\ell_1) \Gamma(\rho_2/\ell_2)}
                {\Gamma(\rho_1/\ell_1 + \rho_2/\ell_2 + k + 1)}
    \Bigr\vert,
\end{equation*}
so that, by the symmetry of the zeros of the Riemann zeta-function, we
have
\[
  \Sigma
  =
  2
  \doublesum_{\substack{\rho_{1} \colon \gamma_{1} > 0 \\ \rho_{2} \colon \gamma_2 > 0}}
    \Bigl\vert
    \frac{\Gamma(\rho_1/\ell_1) \Gamma(\rho_2/\ell_2)}
                {\Gamma(\rho_1/\ell_1 + \rho_2/\ell_2 + k + 1)}
    \Bigr\vert
  +
  2
  \doublesum_{\substack{\rho_{1} \colon \gamma_{1} > 0 \\ \rho_{2} \colon \gamma_2 > 0}}
    \Bigl\vert
    \frac{\Gamma(\rho_1/\ell_1) \Gamma(\overline{\rho}_2/\ell_2)}
                {\Gamma(\rho_1/\ell_1 + \overline{\rho}_2/\ell_2 + k + 1)}
    \Bigr\vert
  =
  2 (\Sigma_1 + \Sigma_2),
\]
say.
It is clear that if both $\Sigma_1$ and $\Sigma_2$ converge, then the
double sum on the right-hand side of \eqref{def-M5}
converges absolutely.
In order to estimate $\Sigma_1$ we choose a large $T$ and let
\begin{align*}
  D_0
  &=
  \{ (\rho_1, \rho_2)\colon (\gamma_1, \gamma_2) \in [0, 2 T]^2 \},
  &
  D_3
  &=
  \{ (\rho_1, \rho_2)\colon  
    \gamma_2 \ge T, \,
    T \le \gamma_1 \le \gamma_2  \}, \\
  D_1
  &=
  \{ (\rho_1, \rho_2)\colon  
    \gamma_1 \ge T, \,
    T \le \gamma_2 \le \gamma_1 \}, 
  &
  D_4
  &=
  \{ (\rho_1, \rho_2)\colon  
    \gamma_2 \ge T, \,
    0 \le \gamma_1 \le T  \}, \\
  D_2
  &=
  \{ (\rho_1, \rho_2)\colon  
    \gamma_1 \ge T, \,
    0 \le \gamma_2 \le T \},
\end{align*}
so that
$\Sigma_1 \le \Sigma_{1,0} + \Sigma_{1,1} + \Sigma_{1,2} + \Sigma_{1,3} +
 \Sigma_{1,4}$,
say, where $\Sigma_{1,j}$ is the sum with $(\rho_1, \rho_2) \in D_j$.
Now, $D_0$ contributes a bounded amount, that depends only on $T$,
and, by symmetry again, $\Sigma_{1,1} = \Sigma_{1,3}$ and
$\Sigma_{1,2} = \Sigma_{1,4}$.
We also recall the inequality \eqref{ineq-zeros}
which is valid for all couples of zeros considered in $\Sigma_1$.
Hence
\begin{align*}
  \Sigma_{1,1}
  &\ll_{\ell_1,\ell_2}
  \doublesum_{\substack{\rho_{1} \colon \gamma_1 \ge T \\ \rho_{2} \colon T \le \gamma_2 \le \gamma_1}}
    \frac{e^{- \pi (\gamma_1/\ell_1 + \gamma_2/\ell_2) / 2} 
    (\frac{\gamma_1}{\ell_1})^{\beta_1/\ell_1 - 1/2}  (\frac{\gamma_2}{\ell_2})^{\beta_2/\ell_2 - 1/2}}
         {e^{- \pi (\gamma_1/\ell_1 + \gamma_2/\ell_2) / 2} (\gamma_1/\ell_1 + \gamma_2/\ell_2)^{\beta_1/\ell_1 + \beta_2/\ell_2 + k + 1/2}} 
         \\
  &
  \ll_{\ell_1,\ell_2}
  \doublesum_{\substack{\rho_{1} \colon \gamma_1 \ge T \\ \rho_{2} \colon T \le \gamma_2 \le \gamma_1}}
    \frac1{\gamma_1^{1/2} \gamma_2^{1/2} (\gamma_1 + \gamma_2)^{k + 1/2}}
         \\
   & 
   \ll_{\ell_1,\ell_2}
  \sum_{\rho_{1} \colon \gamma_1 \ge T}
    \frac1{\gamma_1^{k + 1}}
    \sum_{\rho_{2} \colon T \le \gamma_2 \le \gamma_1}
      \frac1{\gamma_2^{1/2}}
  \ll_{\ell_1,\ell_2}
  \sum_{\rho_{1} \colon \gamma_1 \ge T}
   \frac{\log \gamma_1}{\gamma_1^{k + 1 / 2}}.
\end{align*}
A similar argument proves that
\[
  \Sigma_{1,2}
  \ll_{k,T, \ell_1,\ell_2} 
  \sum_{\rho_{1} \colon \gamma_1 \ge T}
  \frac1{\gamma_1^{k + 1}},
\]
since $\Gamma(\rho_2)$ is uniformly bounded, in terms of $T$, for
$(\rho_1, \rho_2) \in D_2$.
Summing up, we have
\[
  \Sigma_1
  \ll_{k,T, \ell_1,\ell_2} 
  1
  +
  \sum_{\rho_{1} \colon \gamma_1 \ge T}
  \frac{\log \gamma_1}{\gamma_1^{k + 1/2}},
\]
which is convergent provided that $k > 1 / 2$.
In order to estimate $\Sigma_2$ we use a similar argument.
Choose a large $T$ and for $\{ i, j \} = \{1, 2 \}$ set
\begin{align*}
  E_0(i, j)
  &=
  \Bigl\{ (\rho_1, \rho_2) \colon
    \Bigl(\frac{\gamma_i}{\ell_i}, \frac{\gamma_j}{\ell_j} \Bigr) \in [0, 2 T]^2
  \Bigr\}, \\
  E_1(i, j)
  &=
  \Bigl\{ (\rho_1, \rho_2) \colon  
    \frac{\gamma_i}{\ell_i} \ge 2 T, \,
    0 \le \frac{\gamma_j}{\ell_j} \le T \Bigr\}, \\
  E_2(i, j)
  &=
  \Bigl\{ (\rho_1, \rho_2)\colon  
    \frac{\gamma_i}{\ell_i} \ge 2 T, \,
    T \le \frac{\gamma_j}{\ell_j} \le \frac{\gamma_i}{\ell_i} - T \Bigr\}, \\
  E_3(i, j)
  &=
  \Bigl\{ (\rho_1, \rho_2)\colon  
    \frac{\gamma_i}{\ell_i} \ge 2 T, \,
    \frac{\gamma_i}{\ell_i} - T \le \frac{\gamma_j}{\ell_j} \le \frac{\gamma_i}{\ell_i} \Bigr\},
\end{align*}
so that
$\Sigma_2 \le \Sigma_0(1,2) + \Sigma_1(1, 2) + \Sigma_2(1,2) +
 \Sigma_3(1,2) + \Sigma_3(2,1) + \Sigma_2(2,1) + \Sigma_1(2,1)$,
say, where $\Sigma_r(i,j)$ is the sum with $(\rho_1, \rho_2) \in E_r(i,j)$.
Now, $E_0$ contributes a bounded amount, that depends only on $T$,
$\ell_1$ and $\ell_2$.
We remark that similar arguments apply when dealing with
$\Sigma_1(1,2)$ and $\Sigma_1(2,1)$; $\Sigma_2(1,2)$ and $\Sigma_2(2,1)$;
$\Sigma_3(1,2)$ and $\Sigma_3(2,1)$ respectively.
Again we use \eqref{Stirling} as above; hence
\begin{align*}
  \Sigma_2(1,2)
  \ll_{\ell_1,\ell_2}
  \Bigl(
    \doublesum_{\substack{(\rho_1, \rho_2) \in E_2(1,2) \\ \gamma_2 \le \gamma_1^{1/2}}}
    +
    \doublesum_{\substack{(\rho_1, \rho_2) \in E_2(1,2) \\ \gamma_2 > \gamma_1^{1/2}}}
  \Bigr)
    \frac{ (\frac{\gamma_1}{\ell_1})^{\beta_1/\ell_1 - 1/2} 
     (\frac{\gamma_2}{\ell_2})^{\beta_2/\ell_2 - 1/2} e^{-\pi \gamma_2/\ell_2}}
         {(\gamma_1/\ell_1 - \gamma_2/\ell_2)^{\beta_1/\ell_1 + \beta_2/\ell_2 + k + 1/2}}.
\end{align*}
We bound the first sum by a futher subdivision of the zeros $\rho_2$,
treating differently those with $\beta_2 < \ell_2 / 2$ and the other
ones, if any.
The first sum is
\begin{align*}
  &\ll_{\ell_1,\ell_2}
  e^{-\pi T} \,
  \sum_{\gamma_1 \ge 2 T \ell_1}
    \gamma_1^{\beta_1/\ell_1 - 1/2}
    \sum_{\gamma_2 \in [T \ell_2, \gamma_1^{1/2}]}
      \frac{\gamma_2^{\beta_2/\ell_2 - 1/2}}
           {\gamma_1^{\beta_1/\ell_1 + \beta_2/\ell_2 + k + 1/2}} \\
  &\ll_{T,\ell_1,\ell_2}
  \sum_{\gamma_1 \ge 2 T \ell_1}
    \frac1{\gamma_1^{k + 3/2}}
    \Bigl(
      \sum_{\substack{\beta_2 < \ell_2 / 2 \\ \gamma_2 \in [T \ell_2, \gamma_1^{1/2}]}}
      +
      \sum_{\substack{\beta_2 \ge \ell_2 / 2 \\ \gamma_2 \in [T \ell_2, \gamma_1^{1/2}]}}
    \Bigr)
      \Bigl( \frac{\gamma_2}{\gamma_1} \Bigr)^{\beta_2/\ell_2 - 1/2} \\
  &\ll_{T,\ell_1,\ell_2}
  \sum_{\gamma_1 \ge 2 T \ell_1}
    \frac1{\gamma_1^{k + 3/2}}
    \Bigl(
      \sum_{\substack{\beta_2 < \ell_2 / 2 \\ \gamma_2 \in [T \ell_2, \gamma_1^{1/2}]}}
        \Bigl( \frac{\gamma_1}T \Bigr)^{1/2 - \beta_2/\ell_2}
      +
      \gamma_1^{1/2} \log \gamma_1
    \Bigr)
  \ll_{T,\ell_1,\ell_2}
  \sum_{\gamma_1 \ge 2 T \ell_1}
    \frac{\log \gamma_1}{\gamma_1^{k + 1/2}}.
\end{align*}
The rightmost series over zeros plainly converges for $k > 1 / 2$.
The second sum is
\begin{align*}
  &\ll_{T,\ell_1,\ell_2}
  \sum_{\gamma_1 \ge 2 T \ell_1}
    \gamma_1^{\beta_1/\ell_1 - 1/2}
     e^{-\pi \gamma_1^{1/2} / \ell_2} \,
     \sum_{\gamma_2 \in [\gamma_1^{1/2}, (\gamma_1 / \ell_1 - T) \ell_2]}
       \frac{\gamma_2^{\beta_2/\ell_2 - 1/2}}
         {(\gamma_1/\ell_1 - \gamma_2/\ell_2)^{\beta_1/\ell_1 + \beta_2/\ell_2 + k + 1/2}} \\
  &\ll_{T,\ell_1,\ell_2}
  \sum_{\gamma_1 \ge 2 T \ell_1}
    \gamma_1^{\beta_1/\ell_1 - 1/2}
     e^{-\pi \gamma_1^{1/2} / \ell_2} \,
     (\gamma_1 \log \gamma_1)
     T^{-(\beta_1/\ell_1 + k + 1/2)}
     \gamma_1^{1/2},
\end{align*}
which is very small.
The contribution of zeros in $E_1(1,2)$ is treated in a similar fashion,
using the uniform upper bound $\Gamma(\rho_2) \ll_T 1$, and is also
small.
We now deal with $\Sigma_3(1,2)$: we have
\begin{align*}
  \Sigma_3(1,2)
  &\ll_{\ell_1,\ell_2} 
  \doublesum_{(\rho_1, \rho_2) \in E_3}
    e^{- \pi \gamma_1 / (2\ell_1)} \gamma_1^{\frac{\beta_1}{\ell_1} - \frac12}
    e^{- \pi \gamma_2 / (2\ell_2)} \gamma_2^{\frac{\beta_2}{\ell_2}- \frac12}
    \Bigl(
      \min_{\substack{k + 1 \le x \le k + 3 \\ 0 \le t \le T}} \vert \Gamma(x + i t) \vert
    \Bigr)^{-1} \\
  &\ll_{k, T,\ell_1,\ell_2}
  \sum_{\rho_{1} \colon \gamma_1 \ge 2 T \ell_1}
    e^{- \pi \gamma_1/\ell_1} \gamma_1^{\beta_1/\ell_1 + 1/\ell_1}
    \log(\gamma_1 + T),
\end{align*}
provided that $T$ is large enough.
Here we are using Theorem 9.2 of Titchmarsh \cite{Titchmarsh1986} with
$T$ large but fixed.
The series at the extreme right is plainly convergent.

\subsection*{Acknowledgements}
The second Author gratefully acknowledges financial support received
from PRIN-2015 project ``Number Theory and Arithmetic Geometry.''

\end{document}